\pdfoutput=1
\documentclass[final]{siamltex}
\usepackage{amsmath,amssymb,mathabx,graphicx,epsfig,euscript,bm,amsfonts}
\usepackage[algoruled,vlined,linesnumbered]{algorithm2e}

\SetAlFnt{\small}
\SetAlCapFnt{\small}
\SetAlCapNameFnt{\small}
\SetAlCapHSkip{0pt}
\IncMargin{-\parindent}
\SetKwComment{Comment}{\%\ }{}
\DontPrintSemicolon
\title{Green's function-based time stepping for the Kuramoto-Sivashinsky initial-boundary value problem\thanks{This work was supported by the National Sciences and Engineering Research Council of Canada and the Japan Society for Promotion of Science}}
\author{L. van Veen\thanks{Faculty of Science, University of Ontario
    Institute of Technology, 2000 Simcoe St. N., Oshawa, ON L1H 7K4, Canada ({\tt lennaert.vanveen@uoit.ca})}}

\newcommand{\rf}[1]{(\ref{#1})}
\begin{document}
\maketitle

\begin{abstract}
Both theoretical and numerical studies of the Kuramoto-Sivashinsky equation have mostly considered 
periodic boundary conditions. In this setting, the Fourier decomposition of the solution is central
to theoretical ideas, such as renormalization group arguments, as well as to numerical solution,
allowing for the construction of accurate and efficient
time-steppers using standard pseudo-spectral methods. In contrast, fixed boundary
conditions induce boundary layers and necessitate the use of non-uniform grids, usually generated by
orthogonal polynomials. On such bases, numerical differentiation is ill-conditioned and can potentially 
lead to a catastrophic blow-up of round-off error. In this paper, we use ideas recently explored by 
Viswanath ({\sl J. Comput. Phys.} {\bf 251} (2013), pp. 414-431) to completely eliminate numerical 
differentiation and linear solving from the time-stepping algorithm. We use the Green's function-based method to investigate
elements of the Kuramoto-Sivashinsky dynamics over a range of five decades of the viscosity.
\end{abstract}



\pagestyle{myheadings}
\thispagestyle{plain}
\markboth{L. van Veen}{Green's function based time-stepping}

\section{Introduction\label{intro}}

Boundary conditions form an important part of models of continuous, space and time dependent processes. Prescribing
continuous variables like velocity or concentration, or their derivatives, at material boundaries often
leads to boundary layers, i.e. regions with steep gradients of the dependent variables. Periodic boundary
conditions, in contrast, often lead to more spatially homogeneous solutions.
Schemes for numerical simulation must take the peculiar structure of solutions, induced by boundary conditions, into account.

Periodic boundary conditions arise naturally, and exactly, in some geometries, but are also often used when the actual domain
is practically infinite or no more natural choice is available. It is straightforward to list a number of reasons why
periodic boundary conditions are pleasant to deal with.
Firstly, the natural choice of a grid for spatial discretisation is
a regular grid, and the natural choice of a basis to expand the solution in is one consisting of sines and cosines.
The Fast Fourier Transform (FFT) provides an efficient way to switch between grid and spectral representations of a solution.
Each basis function satisfies the boundary conditions and is an eigenfunction of any spatial differential operator.
Exploiting these facts, we can use standard pseudo-spectral methods to simulate the system. 
Examples of pseud-spectral methods for semilinear partial differential equations, such as the Kuramoto-Sivashinsky (KS)
equation, can be found,
for instance, in Trefethen's text book \cite{Trefethen_book}. The essential steps can be summarised as 
\begin{enumerate}
\item Discretize time using an implicit method for linear terms and an explicit method for nonlinear terms.
This leads to a periodic Boundary Value Problem (BVP) to solve for each time step.
\item Formulate the BVP in Fourier space. The linear differential operator is represented by a diagonal matrix, so that
the solution can be written explicitly in terms of the Fourier transform of the nonlinear terms.
\item Use the inverse FFT to find the dependent variables and their derivatives on a regular grid, evaluate the nonlinear
terms there and use the FFT to find their contribution to the BVP.
\end{enumerate}
On a fixed number of $n$ grid points in each spatial dimension, this approach yields a time step that requires $\mbox{O}(n\ln n)$ FLoating point
OPperations (FLOPs) for each spatial dimension and has a bound on the spatial discretization error of the form $\exp(-c n)$,
for come positive constant $c$,
provided that the solution is smooth. 

With fixed boundary conditions, some complications arise. The Fourier bases and regular grids are no longer optimal,
as they leads to spurious oscillations known as Gibb's phenomenon. Instead, theory prescribes the use of orthogonal
polynomials and clustered grids generated by their zeros or extrema. These basis functions are not eigenfunctions of
all spatial differential operators and do not usually satisfy the boundary conditions. The most straightforward way to
address these issues is to adjust the above scheme as follows:
\begin{enumerate}
\item[2a] Formulate the BVP in spectral space. The linear differential operator is represented by a matrix that
is structured (e.g. triangular) but not diagonal.
\item[2b] Delete one linear equation for each boundary condition and replace it by the boundary condition on the polynomial
basis.
\item[3] Use an appropriate (inverse) fast transform to evaluate the nonlinear terms on the clustered grid.
\item[4] Solve the resulting linear problem.
\end{enumerate}
In step 2b, only a small error is incurred for well-resolved solutions. Furthermore, for many polynomial
bases, fast transforms from grid point values to and from expansion coefficients are know. For Chebyshev
bases, for instance, the FFT can be used, as explained by Trefethen \cite{Trefethen_book}, who presents several
explicit examples of this approach.

The most important issue, however, is the fact that the algorithm now requires a solver for a large linear system which is 
ill-conditioned for high-order differentiation and fine grids. For differentiation on Chebyshev bases
the condition number grows as $n^{2 p}$, where $p$ is the order of differentiation. The KS equation
is of fourth order and thus the condition number in a naive application of the Chebyshev spectral algorithm 
does not yield a useful bound on the accumulation of errors. 

These are various ways to improve the conditioning of the linear system. It appears that the most successful approach
combines two elements: preconditioning the linear system by a spectral integration matrix and expanding the highest
derivative of the unknown functions, rather than the unknowns themselves, in a polynomial series. Muite \cite{muite}
compares several subtly different variations of this approach on grids up to $10,000$ points. The preconditioned
system still has a condition number that grows algebraically with the number of grid points and is of no use
when estimating errors. A detailed examination of the linear systems arising in the various forms of spectral
integration by Viswanath \cite{viswanath2} lead to the conclusion that some can yield results near machine accuracy
in spite of the bad conditioning. This good accuracy hinges on the cancellation of errors and the careful implementation
of the boundary conditions.

In the current paper, we follow a different, more radical approach, which was also investigated by Viswanath \cite{Viswanath}.
We solve the BVP analytically using Green's function, thereby eliminating the need for numerical differentiation and
linear solving altogether. Instead, we must construct a proper quadrature, that has an exponentially small error in spite of
the finite differentiability of Green's function.

We show, that a combination of barycentric re-interpolation to sub grids and local Clenshaw-Curtis quadrature leads 
to an accurate time-stepping algorithm that is stable to over resolution and requires 
only moderate computational resources to handle grids with tens of thousands of points. This enables us to generate
numerical solutions to the KS Initial Boundary Value Problem (IBVP) for small values of the viscosity or, equivalently, on very large domains. The goal
of those computations is to generate high quality data on the statistical behaviour of the solutions to the KS equation.
In particular, we would like to find evidence for, or against, a conjecture put forward by Yakhot \cite{Yakhot}, which
states that certain average properties of transient behaviour are described by the exponents of the Kardar-Parisi-Zhang
equivalence class \cite{Kardar}. To the best of our knowledge, no other time-stepping algorithms with spectral accuracy
have been developed for this purpose

\section{The KS IBVP\label{KS_IBVP}}

We consider the following IBVP, originating in the work of Kuramoto \cite{Kuramoto} and 
Sivashinsky \cite{Sivashinsky}:
\begin{align}
u_t+u u_x +u_{xx}+\nu u_{xxxx} &= 0 \label{KSIBVP}\\
u(-1)=l,\ u(1)=r,\ u_{xx}(-1)=u_{xx}(1)&=0 \label{BC1}
\end{align}
Other boundary conditions can of course be considered, fixing for instance $u$ and $u_{x}$ at the
boundaries. In that case, the time-stepping algorithm described here does not change, but Green's
function does, i.e. the analysis presented in Appendix \ref{apA} must be adjusted. 
In the literature, the KS equation is often presented in the scaling 
$\bar{u}=\sqrt{\nu}\,u$, $\bar{x}=(x+1)/\sqrt{\nu}$ and $\bar{t}=t/\nu$, which yields an equation identical
to \rf{KSIBVP} but with $\nu\equiv 1$, considered on $[0,L]=[0,2/\sqrt{\nu}]$. This scaling is
used, for instance, in the numerical work of Wittenberg and Holmes \cite{Wittenberg} and the theoretical
work of Galaktionov {\sl et al.} \cite{Galaktionov}. The former work provides an overview of typical
dynamics generated by the KS equation and a rich reference list.
In the latter work, the existence of a bounded solution 
for any finite time of this IBVP is proven.

Another common guise is the integral formulation
\begin{equation}\label{dhdt}
h_t=\frac{1}{2}(h_x)^2-h_{xx}-\nu h_{xxxx}
\end{equation}
where $-h_x=u$. This form is often used in the physics literature when 
considering the KS equation as a model for interface growth, for instance in Refs. \cite{Yakhot,Cross}.
For our time-stepping algorithm it is more convenient to solve for $u$ on the fixed spatial domain
and consider the viscosity $\nu$ as control parameter. We then typically observe a transition from
equilibrium at large values of $\nu$ to time-periodic motion and finally spatio-temporal chaos as
$\nu$ is decreased. The spatio-temporal chaos exhibits a form of extensivity, as demonstrated in section \ref{scaling}.

\section{Time-stepping based on Green's function\label{time-stepping}}

In this section, we will first describe the Green's function-based algorithm on a high level, before specifying and justifying 
our choice of grids, interpolation and quadrature rules.

It is convenient to consider deviations from a linear profile:
\begin{align}
v_t+v v_x+R\, v+\phi v_x+ R\, \phi+v_{xx}+\nu\, v_{xxxx}&=0  \label{KSIBVPHDBC}\\
v(-1)=v(1)=v_{xx}(-1)=v_{xx}(1)&=0 \label{HDBC}\\
u=v+\phi,\ \phi(x)=l+R (x+1),\,R=\frac{r-l}{2} &
\end{align}
which yields an IBVP with homogeneous Dirichlet boundary conditions. 

\subsection{Reformulation as a linear BVP\label{linearBVP}}

The first step is to turn the IBVP into a linear BVP for every time step. This is done by the use of
an implicit-explicit time discretization. In particular, we use
a Semiimplicit Backward Differentiation Formula (SBDF) \cite{Ascher} to obtain
\begin{equation}\label{SBDF1}
\mathcal{L}v_{k+1}=\sum\limits_{s=1}^{o}\alpha_s v_{k+1-s} + \Delta \sum\limits_{s=1}^{o} \beta_{s} f(v_{k+1-s})
\end{equation}
where $\Delta$ is a constant that depends on the time step size $h$ and on $o$, the order of the SBDF formula. The sub scripts denote the
approximate solution at different times. We have introduced the linear operator $\mathcal{L}$ and the
part of the BVP that is treated explicitly, $f$, according to
\begin{align}
\mathcal{L} &= 1+\Delta R+\Delta \partial_{xx}+\Delta \nu \partial_{xxxx} \label{defL}\\
f(v)&=-v v_x -\phi(x) v_x - R \phi(x) \label{deff}
\end{align}
The definitions of $\Delta$ and the coefficients $\alpha$ and $\beta$ are listed in Table \ref{SBDF_table_1}.
\begin{table}\label{SBDF_table_1}
\begin{center}
\begin{tabular}{l|l|lllll}
      & $\Delta$   &      & $k=1$ & $k=2$ & $k=3$ & $k=4$ \\\hline
$o=1$ &  $h$       & $\alpha_k$ & $1$ & & & \\
      &            & $\beta_k$  & $1$ & & & \\\cline{1-5}
$o=2$ & $2h/3$     & $\alpha_k$ & $4/3$ & $-1/3$ & & \\
      &            & $\beta_k$  & $2$ & $-1$   & & \\\cline{1-6}
$o=3$ & $6h/11$    & $\alpha_k$ & $18/11$ & $-9/11$ & $2/11$ & \\
      &            & $\beta_k$  & $3$     & $-3$    & $1$    & \\\cline{1-7}
$o=4$ & $12h/25$   & $\alpha_k$ & $48/25$ & $-36/25$& $16/25$ & $-3/25$   \\
      &            & $\beta_k$  & $4$     & $-6$    & $4$     & $-1$  \\\cline{1-7}
\end{tabular}
\end{center}
\caption{Table of coefficients of SBDF formula \rf{SBDF1} for orders up to $o=4$.}
\end{table}

The time-discretized equation \rf{SBDF1}, together with boundary conditions \rf{HDBC}, constitutes a linear BVP
since all quantities on the right-hand side are known explicitly. The solution can be written in terms of 
Green's function as
\begin{equation}\label{SBDF2}
v_{k+1}=\sum\limits_{s=1}^{o}\alpha_s G \ast v_{k+1-s} + \Delta \sum\limits_{s=1}^{o} \beta_{s} G \ast f(v_{k+1-s})
\end{equation}
where the star denotes the convolution
$$
(G\ast v)(x)=\int\limits_{y=1}^1 G(x,y) v(y)\,\mbox{d}y
$$
This expression is not suitable for numerical quadrature because the second convolution contains derivatives of
the unknown function. We apply integration by parts to obtain
\begin{multline}\label{SBDF3}
v_{k+1}=\sum\limits_{s=1}^{o} (\alpha_s +\beta_s \Delta R) G\ast  v_{k+1-s} + \Delta \sum\limits_{s=1}^{o}\beta_s DG\ast\left(\frac{1}{2} v_{k+1-s}^2+\phi v_{k+1-s}\right)-\Delta R G\ast \phi\\
\equiv G\ast I_1+DG\ast I_2 +J
\end{multline}
where $DG$ denotes the derivative of Green's function with respect to the variable of integration of the convolution
and, for later convenience, we have introduced shorthand notation for the terms that appear in the convolutions with $G$ and $DG$ and
the constant term.
Now the only derivative remaining can be computed analytically so the problem of the bad conditioning of numerical
differentiation has been eliminated. In its place, we face two new challenges. Firstly, we must compute Green's function and
cast it in a form suitable for numerical evaluation. In Appendix \ref{apA}, appropriate explicit expressions are derived. 
Secondly, we must accurately compute the convolutions in the knowledge that the integrand is only once or twice continuously differentiable. 
For this end, we will use polynomial re-interpolation followed by standard quadrature.

\subsection{Re-interpolation and quadrature\label{Interpolation}}

When using classical, global quadrature
rules, the finite differentiability of the integrands would lead to a fixed rate of convergence, much like a standard finite difference method for solving the
IBVP (\ref{KSIBVP},\ref{BC1}) would have given us.
We can solve this issue by using separate quadratures for the sub domains $[-1,x]$ and $[x,1]$ on which $G(x,y)$ is
smooth. This will require reinterpolation of the integrand onto suitable, non-uniform sub grids.  

Let us denote the global grid by $x_i$, $i=0,\ldots, n$ and the left and right sub grids for the $i^{\rm th}$ global grid point by 
$x_j^{{\textsc l},i}$, $j=0,\ldots,n^{{\textsc l},i}$
and $x_j^{{\textsc r},i}$, $j=0,\ldots,n^{{\textsc r},i}$. In addition, we will denote the approximations of a quantity $a(x)$ on these grids by 
$\bm{a}$, $\bm{a}^{{\textsc l},i}$ and $\bm{a}^{{\textsc r},i}$. Then the convolutions with $G$ and $DG$ are represented by matrix-vector products
\begin{align}
b&=G\ast a \ \ \rightarrow\ \ \bm{b} = M \bm{a} \nonumber\\ 
M_{ij}&=\sum\limits_{p,q=0}^{n^{{\textsc l},i}} C_{p}^{{\textsc l},i}G_{pq}^{{\textsc l},i}B_{qj}^{{\textsc l},i}+\sum\limits_{p,q=0}^{n^{{\textsc r},i}} C_{p}^{{\textsc r},i}G_{pq}^{{\textsc r},i}B_{qj}^{{\textsc r},i}\nonumber\\
b&=DG\ast a \ \ \rightarrow\ \ \bm{b} = N \bm{a} \nonumber\\ 
N_{ij}&=\sum\limits_{p,q=0}^{n^{{\textsc l},i}} C_{p}^{{\textsc l},i}\widebar{G}_{pq}^{{\textsc l},i}B_{qj}^{{\textsc l},i}+\sum\limits_{p,q=0}^{n^{{\textsc r},i}} C_{p}^{{\textsc r},i}\widebar{G}_{pq}^{{\textsc r},i}B_{qj}^{{\textsc r},i}
\end{align}
where
\begin{itemize}
\item $B^{{\textsc l},i}$ is a $(n^{{\textsc l},i}+1)\times (n+1)$ matrix of interpolation from the global grid to the $i^{\rm th}$ left sub grid, and likewise for $B^{{\textsc r},i}$,
\item $G^{{\textsc l},i}$ and $\widebar{G}^{{\textsc l},i}$ are the $(n^{{\textsc l},i}+1)\times (n^{{\textsc l},i}+1)$ diagonal matrices with the values $G(x_i,y)$ and $\partial_y G(x_i,y)$ take on 
the $i^{\rm th}$ left sub grid, i.e.
\begin{alignat}{2}
G_{pq}^{{\textsc l},i}&=\delta_{pq}\,G(x_i,x_p^{{\textsc l},i}) &\qquad \bar{G}_{pq}^{{\textsc l},i}&=\delta_{pq}\,\left.\frac{\partial G(x,y)}{\partial y}\right|_{(x_i,x_p^{{\textsc l},i})}
\end{alignat}
and likewise for $G^{{\textsc r},i}$ and $\widebar{G}^{{\textsc r},i}$ and
\item $C^{{\textsc l},i}$ is a row vector of quadrature coefficients such that
$$
C^{{\textsc l},i} \bm{a}^{{\textsc l},i}\approx \int\limits_{x=-1}^{x_i}a(x,t)\,\mbox{d}x
$$
and likewise for $C^{{\textsc r},i}$.
\end{itemize}

\subsection{Explicit expressions\label{explicit}}

Here, we will give explicit expressions and assess the expected degree of accuracy for
the combination of barycentric interpolation on closed Chebyshev grids and Clenshaw-Curtis quadrature.

We will take both the global grid and the sub grids to consist of Chebyshev points of the second kind, i.e.
\begin{align}\label{Chebyshev_grids}
x_i&=\cos\left(\frac{i \pi}{n}\right)\!,\ i=0,\ldots,n \nonumber\\
x_j^{{\textsc l},i}&=(\psi^{{\textsc l},i}\!\circ\cos)\left(\frac{j \pi}{n^{{\textsc l},i}}\right)\!,\ j=0,\ldots,n^{{\textsc l},i}\ \text{with}\ \psi^{{\textsc l},i}(x)=-1+\frac{1}{2}(x_i+1)(x+1) \nonumber\\
x_j^{{\textsc r},i}&=(\psi^{{\textsc r},i}\!\circ\cos)\left(\frac{j \pi}{n^{{\textsc r},i}}\right)\!,\ j=0,\ldots,n^{{\textsc r},i}\ \text{with}\ \psi^{{\textsc r},i}(x)=\frac{1}{2}(1-x_i)(x-1)+1
\end{align}
An explicit expression for the interpolation matrix can then be found in Berrut and Trefethen \cite{Berrut}. It is given by 
\begin{equation}\label{Interpolation_matrix}
B_{pq}^{{\textsc l},i}=\left(\sum\limits_{k=0}^{n}\frac{w_k (x_p^{{\textsc l},i}-x_q)}{w_q (x_p^{{\textsc l},i}-x_k)}\right)^{-1} \ \text{with}\ w_q=\begin{cases}
1/2 & \text{for}\ q=0 \\ (-1)^q & \text{for}\ 0<q<n \\ (-1)^n/2 &\text{for}\ q=n \end{cases}
\end{equation}
and likewise for $B^{{\textsc r},i}$.
The extremal points of the local grids coincide with global grid points. This leads to a singularity in the formula above,
so that we must separately specify that
\begin{alignat*}{2}
B_{0q}^{{\textsc l},i}&=\delta_{qi} &\qquad B_{0q}^{{\textsc r},i}&=\delta_{q0} \\
B_{n^{{\textsc l},i}q}^{{\textsc l},i}&=\delta_{qn} &\qquad B_{n^{{\textsc r}}q}^{{\textsc r},i}&=\delta_{qi} 
\end{alignat*}

The Clenshaw-Curtis row vector is a product of a row vector of quadrature weights with a matrix representing the discrete cosine
transform:
\begin{alignat}{2}\label{Curtis-Clenshaw_vector}
C_p^{{\textsc l},i}&=(x_i+1)\, a_p \sum\limits_{k=0}^{n^{{\textsc l},i}} c_k F_{kp} &\quad \text{where}\ a_p&=\begin{cases}
1/2 &\text{if}\ p=0,n^{{\textsc l},i} \\ 1 &\text{otherwise} \end{cases} \nonumber\\ 
 & & c_k&=\begin{cases}
1 & \text{if}\ k=0 \\ 2/(1-4 k^2) &\text{if $k$ is even and $k>0$} \\
0 & \text{if $k$ is odd} \end{cases} \nonumber \\
 & & F_{kp}&=\frac{1}{n^{{\textsc l},i}}\cos\left(\frac{k p \pi}{n^{{\textsc l},i}}\right)
\end{alignat}
and likewise for $C^{{\textsc r},i}$.

The number points in the local grids, $n^{{\textsc l},i}$ and $n^{{\textsc r},i}$, should be as least as large as the number
of points in the part of the global grid they span, i.e. $n^{{\textsc l},i}>n+1-i$ and $n^{{\textsc r},i}>i+1$. Under this 
condition, the error of interpolation from the global to the local grids is negligible as compared the the error
of interpolation of the solution onto the global grid. Importantly, barycentric interpolation is stable under
over resolution. If we set $n^{{\textsc l},i}=n^{{\textsc r},i}=n+1$ in the implementation described below, the results
remain the same, at least up to machine precision.

\section{Implementation\label{implementation}}

The task of time-stepping the IBVP can now be split up into two steps. In the first step, we compute 
the quadrature matrices $M$ and $N$ for given $\nu$, $\Delta$ and $n$ and a specific choice of the number of points
in the local Chebyshev grids. In the second step we iterate a SBDF by simply computing matrix-vector products.
Both tasks can easily be performed in parallel by distributing the points on the global grid, and
the corresponding rows of the quadrature matrices, over processes. Here, we describe an MPI-based implementation,
which has the advantage that the memory requirements for storing the quadrature matrices can be reduced.
In the following, we label the processes $p=0,\ldots,P-1$, and process number $p$ will compute the solution
on $n^p$ points with indices $i_{\rm s}^{p}$ through $i_{\rm e}^p$. We further assume that $P\ll n$ so that the
work load can be almost evenly balanced, i.e. $n^p\approx n/P$ for $p=0,\ldots,P-1$.
 
Algorithm \ref{Alg:quadrature} shows a pseudo-code for the construction of the quadrature matrices.
The most costly steps in the main loop over global grid points are number 4, the computation of the
Clenshaw-Curtis quadrature weights, and number 6, the computation of the interpolation matrix.
Assuming that the total number of local grid points is of the same order as the number of global grid points,
these steps have a FLOP count of $\mbox{O}(n^2)$, bringing the total for this algorithm to $\mbox{O}(n^3)$.
In the parallel implementation each processor will thus handle $\mbox{O}(n^3/P)$. 

One more remark about the quadrature algorithm must be made. If the order of the global grid is large, 
it may happen that it coincides with points on local grids up to machine accuracy. When using standard double
precision arithmetic, for instance, this start to happen for global grid orders upward from 20,000.
This leads to divisions by zero when evaluating \rf{Interpolation_matrix}. This can simply be circumvented
by detecting overlap up to finite precision of the grids and replacing each row of $B$ that holds
the coefficients of interpolation onto a overlapping local grid point by a row of zero elements and a
single element equal to unity, just like for the extremal points of the local grid. It is a remarkable fact,
explained in detail by Higham \cite{Higham}, that the evaluation of the elements of $B$ is otherwise
stable, in spite of the small denominators. In our test of the accuracy of the numerical quadrature
in section \ref{tests}, we replace rows corresponding to local grid points that are closer to global
grid points than $10^{-13}$ and the resulting computation is stable up to a global grid order of at least
74,000.
\IncMargin{1.5em}
\begin{algorithm}[t]
\KwIn
{
viscosity $\nu$, time-stepping parameter $\Delta$, boundary values $l$ and $r$ and global and local Chebyshev grid orders $n$,
$n^{{\rm L},i}$, $n^{{\rm R},i}$, $i=0,\ldots,n$.
}
\KwOut
{
$(n+1)\times (n+1)$ quadrature matrices $M$ and $N$ and $(n+1)$ vector $J$, distributed by rows over $P$ processes.
}
\For(\Comment*[f]{loop over rows stored on processor $p$})
{
$i=i_{\rm s}^p\ldots i_{\rm e}^p$
}
{
Allocate $a_j^{{\textsc l}},\,b_j^{{\textsc l}},\ j=0,\ldots,n^{{\textsc l},i}$, $a_j^{{\textsc r}},\,b_j^{{\textsc r}},\ j=0,\ldots,n^{{\textsc r},i}$.\;%
Set $x_j^{{\textsc l},i},\ j=0,\ldots,n^{{\textsc l},i}$ and $x_j^{{\textsc r},i},\ j=0,\ldots,n^{{\textsc r},i}$
according to \rf{Chebyshev_grids}.\;%
Compute the Clenshaw-Curtis row vectors according to \rf{Curtis-Clenshaw_vector} and set
\begin{alignat*}{4}
a^{{\textsc l}}&\leftarrow C^{{\textsc l},i} &\qquad b^{{\textsc l}}&\leftarrow C^{{\textsc l},i} &\qquad 
a^{{\textsc r}}&\leftarrow C^{{\textsc r},i} &\qquad b^{{\textsc r}}&\leftarrow C^{{\textsc r},i}
\end{alignat*}\;\vspace{-10pt}
Multiply by Green's function or its derivative, i.e. set
\begin{alignat*}{4}
a^{{\textsc l}} &\leftarrow \sum\limits_{p=0}^{n^{{\textsc l},i}}a_p^{{\textsc l}}G^{{\textsc l},i}_{pq} &\quad
b^{{\textsc l}} &\leftarrow \sum\limits_{p=0}^{n^{{\textsc l},i}}b_p^{{\textsc l}}\widebar{G}^{{\textsc l},i}_{pq} &\quad
a^{{\textsc r}} &\leftarrow \sum\limits_{p=0}^{n^{{\textsc r},i}}a_p^{{\textsc r}}G^{{\textsc r},i}_{pq} &\quad
b^{{\textsc r}} &\leftarrow \sum\limits_{p=0}^{n^{{\textsc r},i}}b_p^{{\textsc r}}\widebar{G}^{{\textsc r},i}_{pq}
\end{alignat*}\;\vspace{-10pt}
Compute the interpolation matrices according to \rf{Interpolation_matrix} and set
\begin{alignat*}{2}
M_{ij}&\leftarrow \sum\limits_{p=0}^{n^{{\textsc l},i}} a_p^{{\textsc l}} B^{{\textsc l},i}_{pj}+
\sum\limits_{p=0}^{n^{{\textsc r},i}} a_p^{{\textsc r}} B^{{\textsc r},i}_{pj} &\qquad 
N_{ij}&\leftarrow \sum\limits_{p=0}^{n^{{\textsc l},i}} b_p^{{\textsc l}} B^{{\textsc l},i}_{pj}+
\sum\limits_{p=0}^{n^{{\textsc r},i}} b_p^{{\textsc r}} B^{{\textsc r},i}_{pj}\quad\text{for}\ j=0,\ldots,n
\end{alignat*}\;
De-allocate $a_j^{{\textsc l}},\,b_j^{{\textsc l}},\,a_j^{{\textsc r}},\,b_j^{{\textsc r}}$.\;
Compute $J_i=-\Delta R \sum_{j=0}^n M_{ij}\phi(x_j)$.\;
}
\Return{$M^p$, $N^p$ and $J^p$, i.e. rows $i_{\rm s}^p$ through $i_{\rm e}^p$ of the quadrature matrices and the constant term.}
\caption{Computation of the quadrature matrices and constant term} \label{Alg:quadrature}
\end{algorithm}
\DecMargin{1.5em}

Algorithm \ref{Alg:time_stepping} describes the time-stepping. The loop over time steps includes
the computation of the integrands in equation \rf{SBDF3}, which takes $\mbox{O}(n^p)$ FLOPS on
processor $p$, an MPI all-to-all communication of $\mbox{O}(n)$ elements and a matrix vector
product with the elementary FLOP count of $\mbox{O}(n^p n)$. If the MPI routine
takes a similar amount of time to complete as $\mbox{O}(n)$ FLOPS, it is obvious that the communication time
will be negligible as long as $n^p\gg 1$, as we have assumed.

\IncMargin{1.5em}
\begin{algorithm}[t]
\KwIn
{time step size $h$, number of inner and outer iterations $N_{\rm i}$, $N_{\rm o}$, $o$ initial points $v_{-o+1},\ldots v_{0}$,
boundary values $l$ and $r$.\\\mbox{\hspace{33pt}}$M^p$, $N^p$ and $J^p$ are stored on processor $p$.}
\For{$j=1,\ldots,N_{\rm o}$}
{
  MPI scatter $v_{k-o+i}\rightarrow v_{k-o+i}^0\ldots v_{k-o+i}^{p-1}$ for $i=1,\dots,o$\Comment*[f]{root to all}\;
  \For{$k=1,\ldots,N_{\rm i}$}
      {
        Compute the integrands $I_1^p$ and $I_2^p$ of \rf{SBDF3}.\;
        MPI gather the integrands: $I_{1,2}^0,\ldots,I_{1,2}^{p-1}\rightarrow I_{1,2}$.\Comment*[f]{all to all}\;
        Set
        \begin{equation*}
          v^p_{k+1}=M^p I_1+N^p I_2+J^p
        \end{equation*}\;
      }
      MPI gather $v_{k}^0,\ldots v_{k}^{p-1}\rightarrow v_k$.\Comment*[f]{all to root}\;
      Root: output $v_k+\phi$.\;
}
\caption{Green's function based time-stepping} \label{Alg:time_stepping}
\end{algorithm}
\DecMargin{1.5em}

\begin{figure}[t]\label{fig:wtime}
\begin{center}
\includegraphics[width=0.49\textwidth]{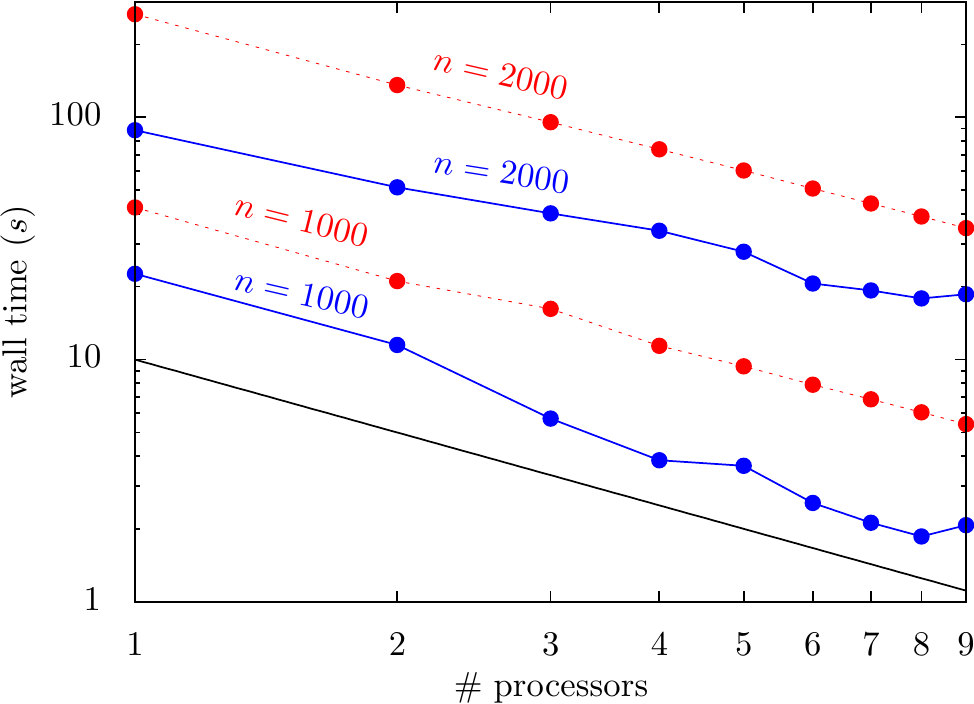}\quad
\includegraphics[width=0.47\textwidth]{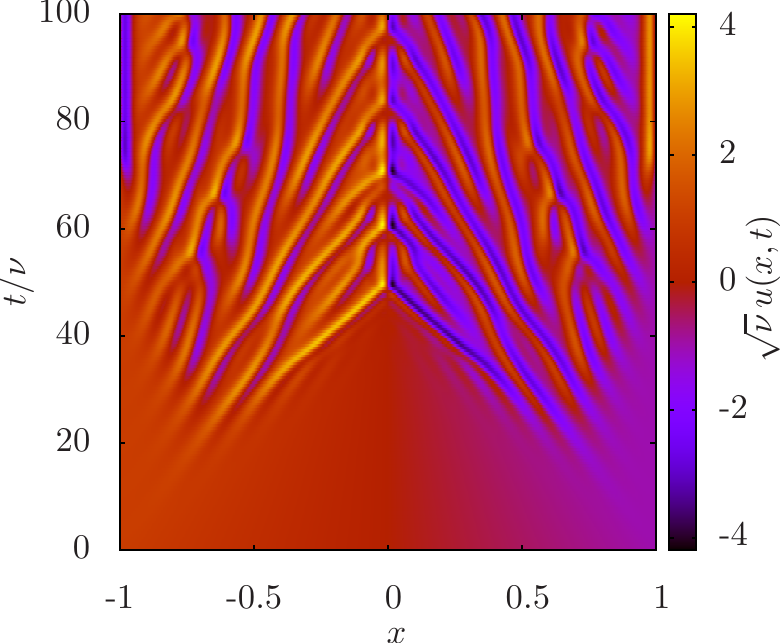}
\end{center}
\caption{Left: wall time for the computation of the quadrature matrices according to Algorithm \ref{Alg:quadrature} (red) and
for time stepping according to Algorithm \ref{Alg:time_stepping} (blue) using the first order SBDF. Shown is the wall time averaged over 20 trial for 4000 time steps
with grid orders $n=1000$ and $n=2000$. The black line indicates linear scaling. Right: density plot of the solution computed in this test. The viscosity is set to
$\nu=2\times 10^{-4}$, the step size to $h=10^{-5}$ and the initial condition is $u_0=-x/\sqrt{\nu}$.}
\end{figure}

Figure \ref{fig:wtime} shows the wall time taken for building the quadrature matrices and time-stepping $4000$ times
for grid sizes $1000$ and $2000$. 
This test was run on a cluster computer, on a node with 24 AMD Opteron 2.2GHz processors, 32Gb of RAM memory,
512Kb of cache and a QDR InfiniBand connection.
Clearly, Algorithm \ref{Alg:quadrature} scales linearly as predicted by the FLOP count since it does not involve
any communication. Time-stepping according to Algorithm \ref{Alg:time_stepping}, on the other hand, shows approximately
linear speed-up up to 9 processors only. This results depends on architecture, on a machine with two quad-core
CPUs, for instance, the speed-up saturates at 2 processors. There are two possible reasons for this limitation. Firstly,
the execution of the all-to-all communication depends on the details of the hardware and the MPI implementation
and is highly sensitive to concurrent use by other processes of cache memory. Secondly, the all-to-all communication
is blocking and the execution time for the partial matrix-vector product may vary over processes.

Lastly, we turn to the nontrivial question of initializing the time-stepping. For SBDF orders two and up, we need
the solution at previous time instants. We propose four possible solutions:
\begin{enumerate}
\item Using the first order SBDF formula with a small time step. This is a commonly used method to seed backward differentiation
formulae, and is employed by Ascher {\sl et al.} \cite{Ascher}. For the Green's function based time-stepper, however, this
method has its limitation. If the time step is taken very small for fixed viscosity, Green's function approaches a delta
distribution as its length scale, $b^{-1}$  (see \rf{auxiliaries}) goes to zero. If this scale becomes comparable to the
spacing of the local grid near its end points, of order $n^{-2}$, an instability can occur. 
\item Richardson extrapolation from lower order. For instance, we can seed the second order SBDF by approximating
$u(x,h)$ first by two first order steps of size $h/2$, giving $u_{h/2}(x,h)$, then by a single step, giving $u_h(x,h)$
and finally setting $u(x,h)\approx 2 u_{h/2}(x,h)-u_h(x,h)$. Similar expressions can readily be found for higher order seeding.
This method has two disadvantages. Firstly, it relies on the cancellation of error terms , which is likely unstable for
the higher order versions, so that round-off error is introduced. Secondly, for each step with a different value of $\Delta$,
we need to repeat Algorithm \ref{Alg:quadrature}, which has order $\mbox{O}(n^3)$ complexity. Therefore, this method is only practical
for small grids and low SBDF order.
\item Using a known exact solution to the KS equation. An exact solution can be found, for instance, in Parkes and Duffy \cite{Parkes}.
We use it in section \ref{tests} to test the accuracy of the SBDF formulae.  Strictly speaking, this is not a solution
to the IBVP, but for small enough viscosity the boundary conditions are satisfied far beyond machine accuracy.
\item Growing a solution from a small perturbation to the zero solution. If we compose the perturbation out of
eigenmodes of the linear part of \rf{KSIBVPHDBC}, we can compute the solution backward in time under the assumption 
that the nonlinear term in negligible. The disadvantage of this method is that it requires a long time integration
for the perturbations to grow to finite size. This strategy is used in our computation of the finite-size effects in
section \ref{scaling}.
\end{enumerate}

\section{Convergence and stability tests\label{tests}}

We present two tests to evaluate the accuracy of the Green's function based time stepping. The first test demonstrates the
exponential convergence of the quadrature computed in each time step. We used the test function
\begin{equation}\label{testfunc}
\xi(x,k)=\frac{1}{1+\sin(\pi k x)^2}-\frac{1}{2}\cos(2\pi k x)-\frac{1}{2}
\end{equation}
which satisfies the homogeneous Dirichlet boundary conditions \rf{HDBC} and describes $2k$ oscillations on the
domain. A similar function was used by Trefethen to demonstrate the convergence of standard spectral methods 
for functions that can be analytically continued in a neighbourhood of $[-1,1]$ in the complex plane \cite{Trefethen_book}.
The continuation of our test function has poles at $\pm i \ln(1+\sqrt{2})/(k \pi)$ that determine an upper bound
for the error of Lagrange interpolation on the global grid \cite{Berrut}, as follows:
$$
\max_{x\in [-1,1]}|\xi(x,k)-p_n(x)|\leq C \exp\left(-\frac{\alpha}{k}n\right)
$$
where $C$ is a positive constant and $p_n$ the interpolant of order $n$. For $k\ge 6$, we find that
$\alpha\approx 0.28$. In fact, $\alpha=\ln(1+\sqrt{2})/\pi$ up to corrections of order $\mbox{O}(1/k^3)$.
This interpolation error sets, in turn, an upper bound for the error in the Clenshaw-Curtis quadrature, the difference
being a constant factor \cite{Trefethen}.
\begin{figure}[t]\label{fig:quad_err}
\begin{center}
\includegraphics[width=0.6\textwidth]{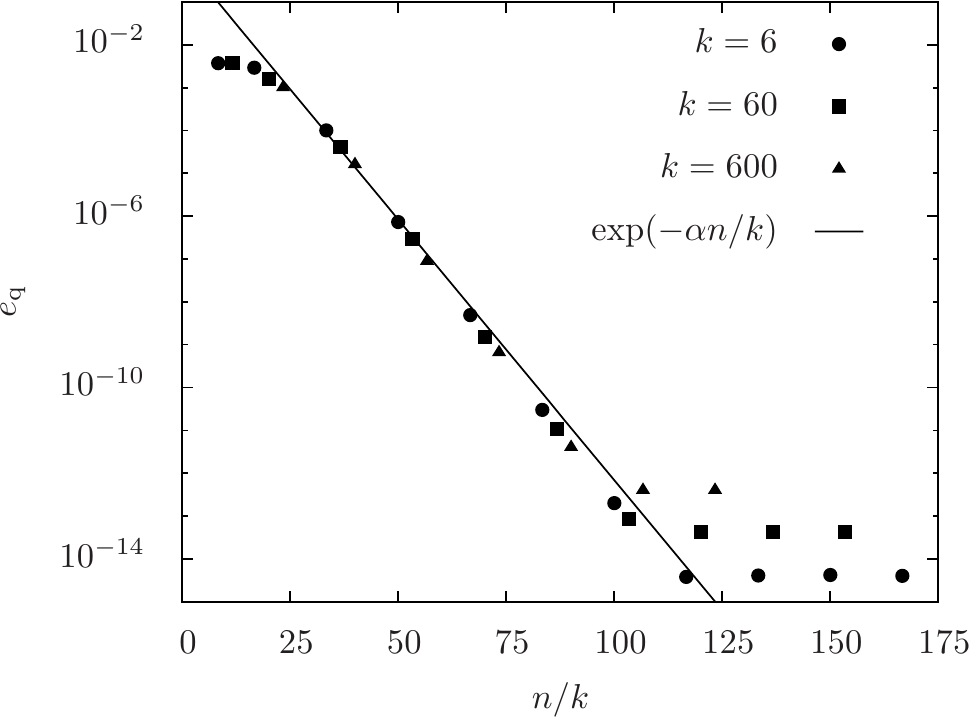}
\end{center}
\caption{Error of the numerical approximation of the convolution with Green's function. Test function \rf{testfunc} 
is considered with $12$, $120$ and $1200$ oscillations in $[-1,1]$. The solid line denotes the theoretically expected
error of the interpolation and Clenshaw-Curtis quadrature.}
\end{figure}

Let $\bm{\xi}$ represent the test function 
evaluated on the global grid, and let $\bm{\omega}$ represent $\mathcal{L}\xi$ evaluated on the global grid. Then we measure
the quadrature error
$$
e_{\rm q}=\|\bm{\xi}-M \bm{\omega}\|_{\infty}
$$
as a function of the grid order $n$ for fixed $\nu$, $h$ and $k$. In the first test, we set $\nu_1=\nu=2\times 10^{-4}$,
$h=\Delta=10^{-5}$ and $k=6$ to generate a function somewhat similar to the final state of the simulation shown in
figure \ref{fig:wtime}(right). In the second test, we set $\nu_2=\nu=2\times 10^{-6}$,
$h=\Delta=10^{-7}$ and $k=60$ and in the third $\nu_3=\nu=2\times 10^{-8}$,
$h=\Delta=10^{-9}$ and $k=600$. The rationale for this choice of parameter values is that we expect the typical
spatial scale of variation of the solutions to decrease as $\sqrt{\nu}$ for small viscosity, as demonstrated in 
section \ref{scaling}.

The three data sets collapse onto a single curve if we plot $e_{\rm q}$ as a function of
$n/k$. Along this curve, the quadrature error decreases approximately as $\exp(-\alpha n/k)$ as predicted. This results
indicates that the main error introduced by the numerical evaluation of the convolutions is that of
polynomial interpolation. We can make two further observations. Firstly, the method is stable to over resolution, as
the error does not increase beyond the minimum around $n/k=120$.  Secondly, the minimal error appears to be determined
by a build-up of round-off error in the matrix-vector product. Every increase in the number of grid points, and thereby
vector elements, by a factor of ten leads to the same increase in the minimal error for $n/k\gtrapprox 120$.

The second test illustrates the order of convergence of the SBDF formulae. For this end, we use the known exact solution
to the KS equation mentioned above, which is given by
\begin{equation}\label{KS_exact}
w(x,t)=c+\frac{15}{19}\frac{\sqrt{11}}{\sqrt{19\nu}} \left[11  \tanh^3\left(\frac{q}{\sqrt{\nu}} (x-c t-x_0)\right)- 9 
\tanh\left(\frac{q}{\sqrt{\nu}} (x-c t-x_0)\right) \right]
\end{equation}
Here, $c$ is the speed of a solitary wave connecting two
constant solutions, $x_0$ is its initial position and $q=\sqrt{11/76}$ is constant. Of course, $w$ is is not an exact solution to 
the initial boundary value problem, but it approaches its left and right limit values at a rate of $\exp(-2q d/\sqrt{\nu})$
a distance $d$ away from the soliton. We have set $\nu=5\times 10^{-5}$, $c=1000$ and $x_0=-0.2$ such that the variation of $w$ at the boundaries, 
and the the magnitude of its second derivative there, are below $10^{-25}$ for $0\leq t/\nu \leq 8$. A similar test was used by Anders \cite{Anders}
{\sl et al.}, but their
viscosity is two orders of magnitude larger and, consequently, they were forced to consider time-dependent boundary conditions.

\begin{figure}[t]\label{fig:step_test}
\begin{center}
\includegraphics[width=0.59\textwidth]{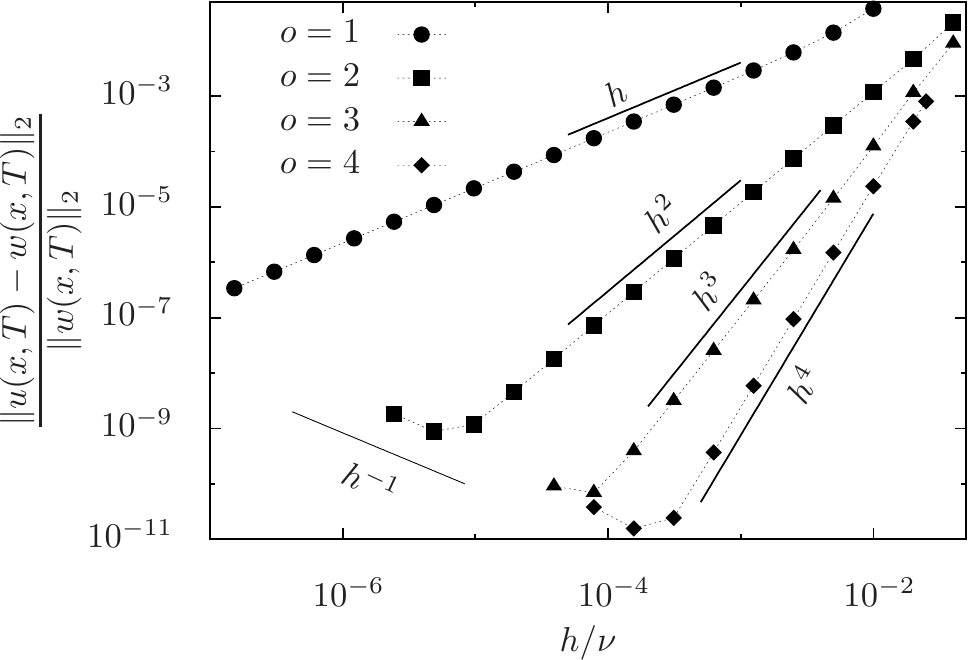}
\includegraphics[width=0.4\textwidth]{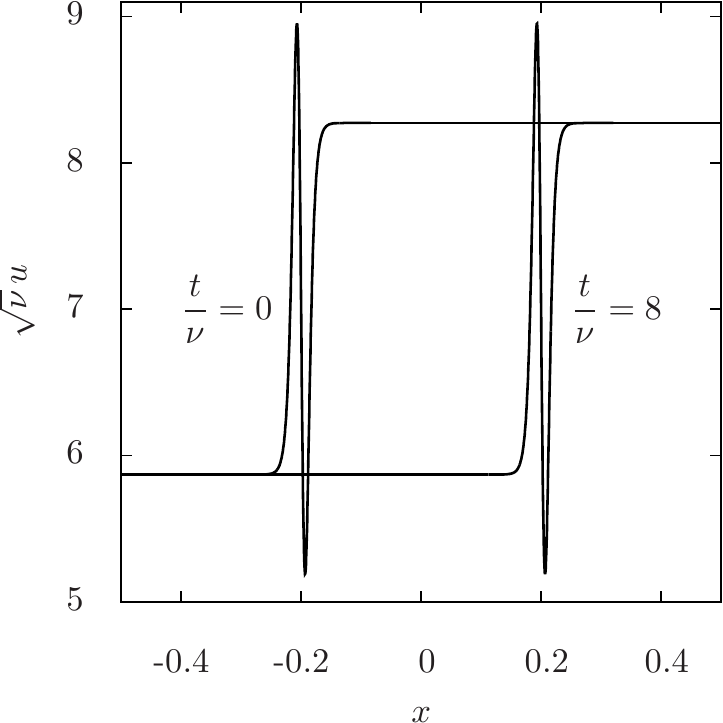}
\end{center}
\caption{Left: relative error in $u$ after time-stepping over $t/\nu=8$ for SBDF formulae of order up to four. The expected
error is indicated by a solid line for each order. Right: initial and final solution in this test.}
\end{figure}
Figure \ref{fig:step_test}(left) shows the error of time-stepping with SBDF formulae up to order four. In these tests, the Chebyshev grid
order was fixed to $n=2000$. For our choice of parameters, the exact solution has a complex singularity close to that of test 
function \ref{testfunc} with $k=6$, and consequently we do not expect the interpolation error to play a role. The smallest error
achieved is instead determined by the accumulation of round-off error, inversely proportional to the step size $h$, as indicated
with a solid line. The initial and final condition in this test are shown in figure \ref{fig:step_test}(right), which shows only
the centre of the domain.

\begin{figure}[t]\label{fig:stab_test}
\begin{center}
\includegraphics[width=0.48\textwidth]{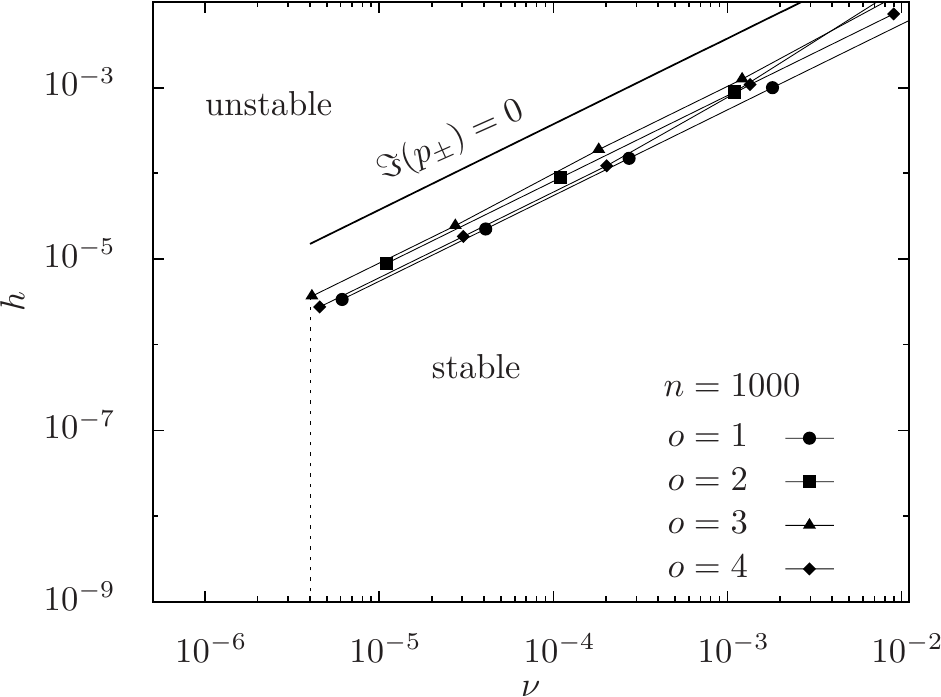}\quad
\includegraphics[width=0.48\textwidth]{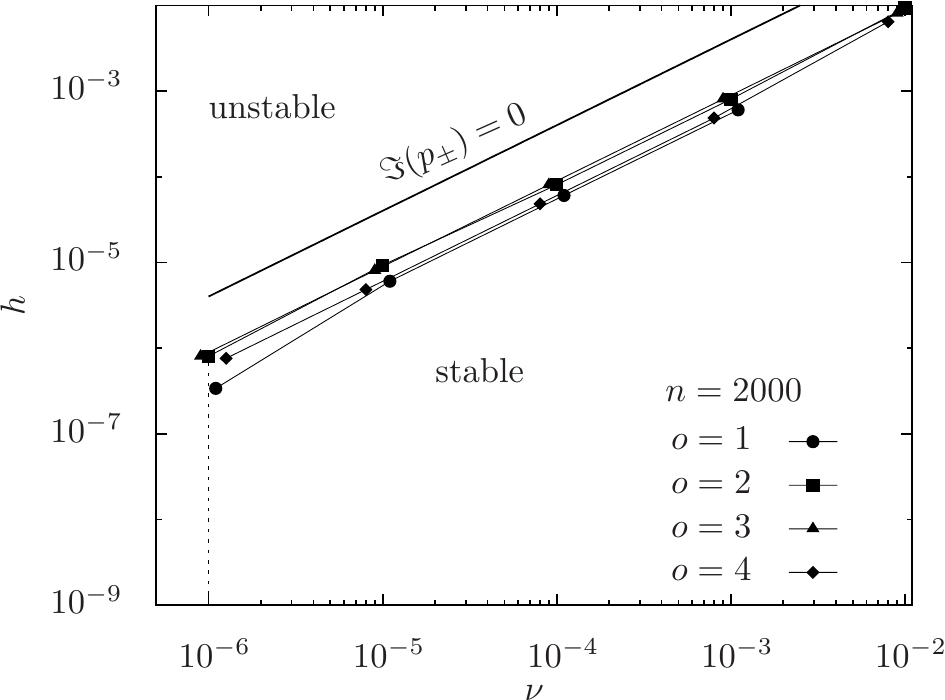}
\end{center}
\caption{Stability of the time-stepping algorithm for SBDF formulae of order up to four. The algorithm is stable below the
solid lines connecting the symbols for each order and to the right of the dashed line. Left: grid order $n=1000$, right $n=2000$.
The solid line at the top of the figures denotes the transition from imaginary to real-valued roots \ref{roots} of the eigenvalues of
the linear operator. The dashed line indicates the point where $\sqrt{\nu}=2/n$, i.e. the largest grid spacing is roughly equal to the
typical length scale of the solution. See text for initial condition.}
\end{figure}

Finally, we consider the stability of the time-stepping method. Figure \ref{fig:stab_test} illustrated the stability for
grid orders $n=1000$ (left) and $n=2000$ (right). The solution was initialized by a linear combination of eigen modes
of the linear operator with small, randomly chosen amplitudes. A combination of viscosity and time step size was labeled
stable if the integration proceeded up to $t/\nu=150$ without blow-up. There are two boundaries to the stable regime.
One lies close to the transition from imaginary to real roots of the eigenvalues of the linear operator given in \ref{roots}.
If these roots are real-valued, Green's function exhibits global oscillations, meaning that updates to the solution become
dependent over an arbitrarily long distance in a single time step, which renders the step unstable. The other occurs
when the typical length scale of the solution, expected to scale as $\sqrt{\nu}$, equals the maximal grid spacing, which is
fixed in these experiments.

\section{Example computations: finite-size effects\label{scaling}}

\begin{figure}[h]\label{fig:contours}
\begin{center}
\includegraphics[width=\textwidth]{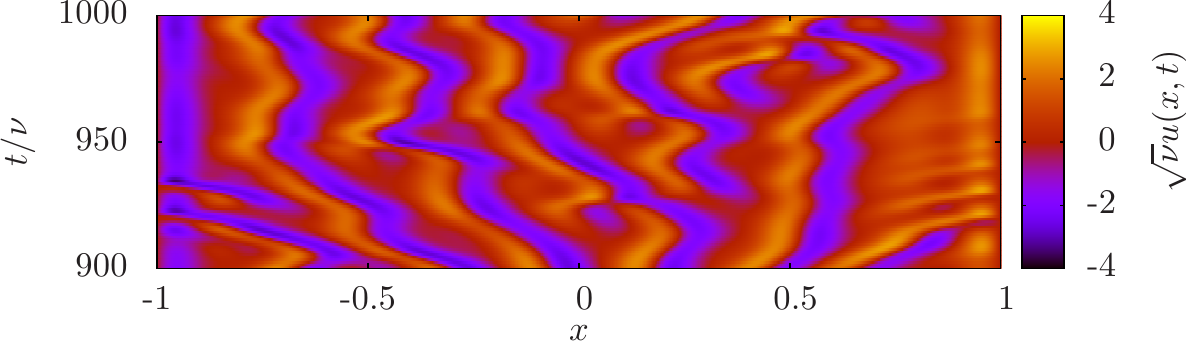}\\
\includegraphics[width=\textwidth]{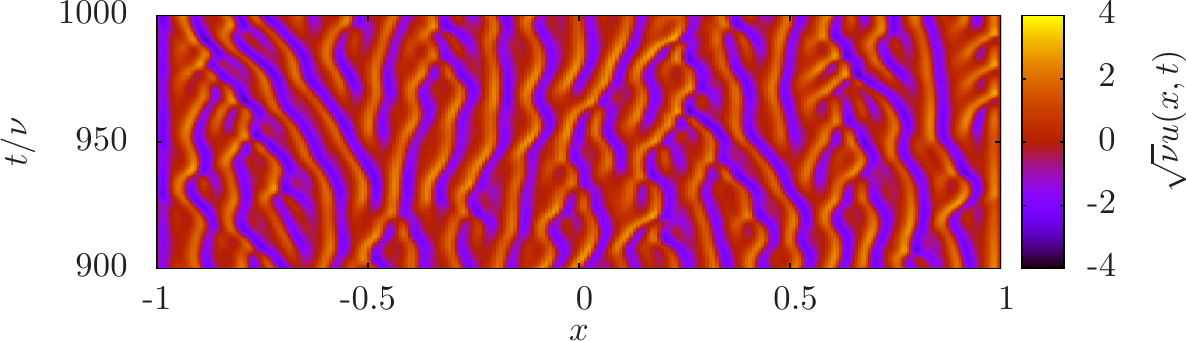}\\
\includegraphics[width=\textwidth]{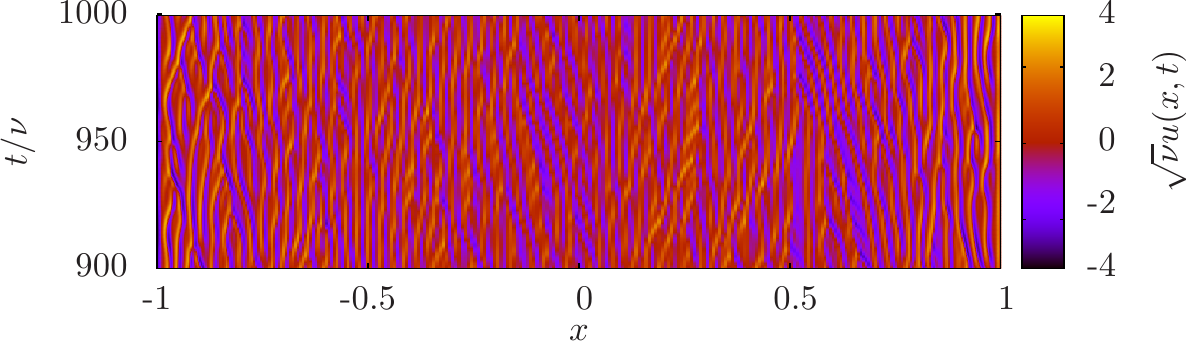}\\
\includegraphics[width=\textwidth]{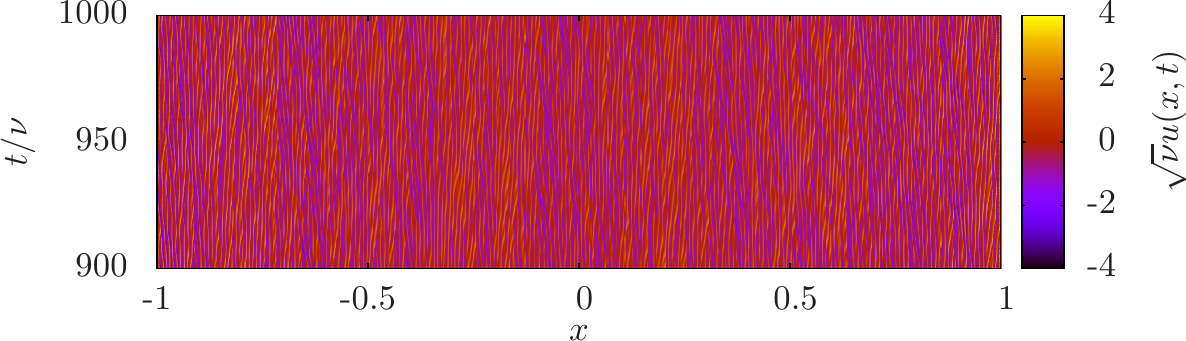}\\
\includegraphics[width=\textwidth]{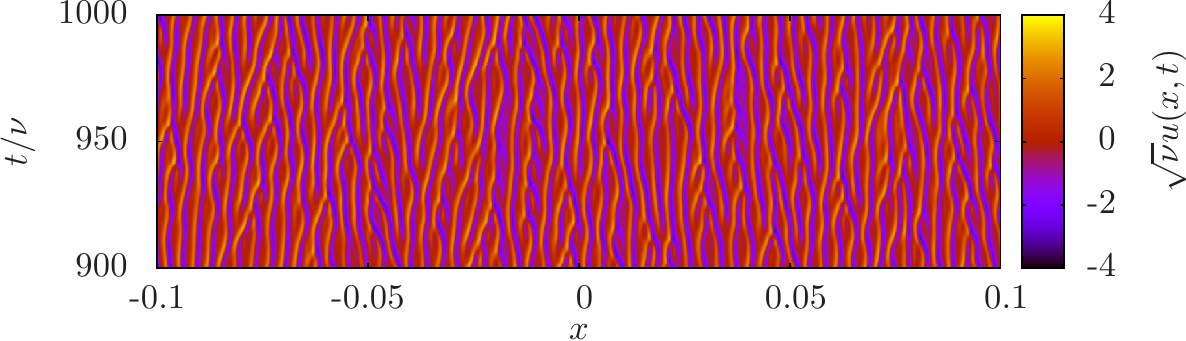}
\end{center}
\caption{Fragment of the time series of $u$ for five orders of magnitude of the viscosity. From top to bottom: $\nu=10^{-3}$ ($n=1000$), $\nu=10^{-4}$ ($n=3000$), 
$\nu=10^{-5}$ ($n=10000$), $\nu=10^{-6}$ ($n=24000$) and $\nu=10^{-7}$ ($n=60000$). Shown are the contours of $\bar{u}=\sqrt{\nu}u$ for $900\leq t/\nu \leq 1000$.
The initial condition in each case is a small random perturbation of the zero solution, the SBDF order is $o=4$ and the time step is $h/\nu=5\times 10^{-4}$.
For $\nu=10^{-7}$, an enlargement of one tenth of the domain is shown so that the qualitative dynamics can be compared to that at $\nu=10^{-5}$.}
\end{figure}

To demonstrate the power of the time-stepping method based on Green's function and Clenshaw-Curtis quadrature, we generated
time series, seeded with random initial conditions of small amplitude, for a range of five orders of magnitude of the viscosity. Each time
series extends up to $t/\nu=2000$ and employed the SBDF formula of fourth order. 
After a transient time of about $t/\nu=150$, the amplitude 
saturates and the dynamics is highly nonlinear. A fragment of each time series is shown in figure \ref{fig:contours}. For the
smallest viscosity, $\nu=10^{-7}$, we have enlarged one tenth of the domain. The fact that the dynamics look qualitatively the same
as that for $\nu=10^{-5}$ is indicative of scaling behaviour, i.e. for small enough viscosity the solutions look similar in the
scaled variables $\bar{u}$, $\bar{x}$ and $\bar{t}$ introduced in section \ref{KS_IBVP}.

\begin{figure}[t]\label{fig:BL}
\begin{center}
\includegraphics[width=0.6\textwidth]{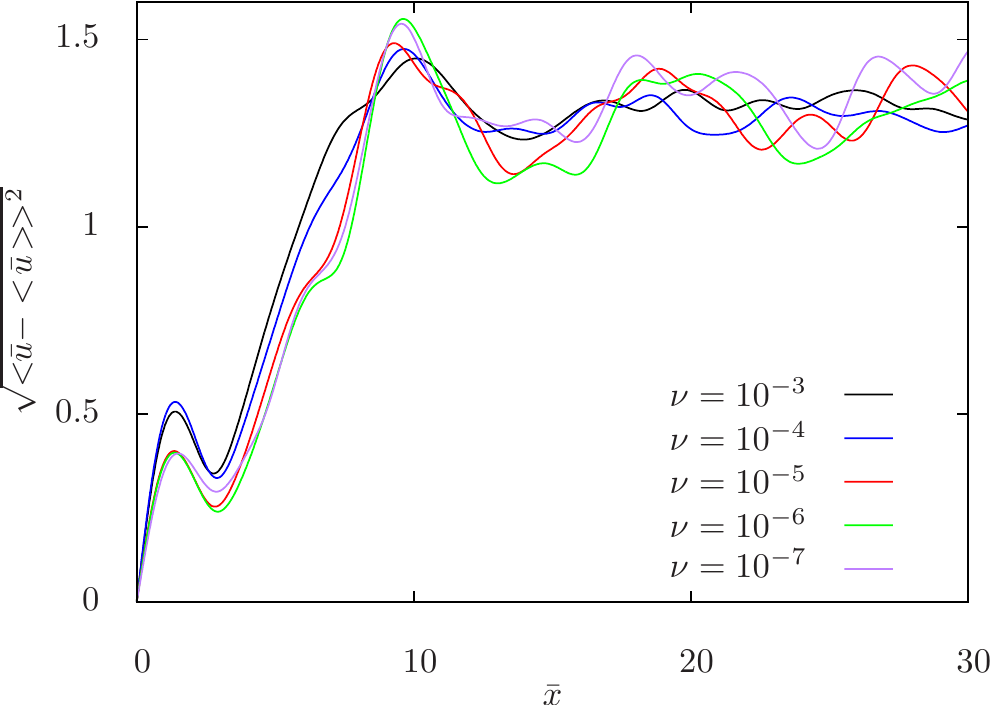}\\
\end{center}
\caption{Root-mean-square of the departure of the solution from the mean profile near the boundary. The brackets denote averaging from $t/\nu=200$ to $t/\nu=2000$
and the scaled variables are $\bar{u}=\sqrt{\nu}u$ and $\bar{x}=(x+1)/\sqrt{\nu}$.}
\end{figure}

Figure \ref{fig:BL} shows the deviation from the time-mean solution near the left boundary in scaled variables to illustrate the
dependence of the boundary layer thickness on viscosity. The fact that the curves approximately
overlap indicates that the thickness of the boundary layer, in which the boundary conditions strongly influence the 
dynamics, is about $12\sqrt{\nu}$ and is constant in scaled variables on the domain $[0,L]$.

\section{Discussion\label{discussion}}

We have demonstrated, that the Green's function based method, in conjunction with the SBDF formulae, barycentric
interpolation and Clenshaw-Curtis quadrature, is accurate, stable and reasonably fast for values of the viscosity
as small as $10^{-7}$. To the best of our knowledge, no other time-stepping method for the KSIBVP has been tested
for a viscosity this small or, equivalently, a domain this large. 

Early studies of the KSIBVP, mostly in the physics literature, used finite-difference discretizations. 
Typical examples of such work are Manneville \cite{Manneville} and Sakaguchi \cite{Sakaguchi}. 
Unfortunately, in this and similar work little description is usually given of the numerical methods, their
accuracy and stability -- Manneville's work being a notable exception. There is no evidence that finite-difference-based
results with a viscosity as low as $10^{-6}$ are reliable even for statistical analysis or curve-fitting exercise like that in
the work cited above.

In later work, various global and piecewise collocation methods were tested on the KSIBVP. Often, validation was only done
for smooth, viscous solutions, like in Khater \& Temsah \cite{Khater}, who used spectral integration on a Chebyshev polynomial
basis. Fornberg \cite{Fornberg} applied a Chebyshev pseudospectral method, implementing the boundary conditions in real
space, and computed chaotic solutions at $\nu\approx 4\times 10^{-4}$.
Piecewise collocation was tested for a viscosity of order $\mbox{O}(10^{-5})$ by Xu \& Shu \cite{Xu} and Anders {\sl et al.}
\cite{Anders}. Based on careful testing and comparison to a priori error estimates, they could conclude that the spatio temporal
chaos they observe numerically is a genuine property of the KS equation rather than a numerical artifact, but it remains unclear if their methods are
suitable for time-stepping at smaller viscosity. 

The limiting factor of the Green's function based method described here is the memory requirement, as it requires storing two $(n+1)\times (n+1)$ 
matrices -- albeit distributes over processors -- and the limited scaling of the parallelization. One possible solution is
to switch to piecewise Chebyshev grids to avoid excessive clustering of grid points near the boundaries. This will require 
a marginally more complicated procedure to compute the quadrature matrices and is work in progress.

\section*{Acknowledgements}
I would like to thank the Faculty of Engineering Science of Osaka University and the Japan Society for Promotion of Science for making possible
the sabbatical visit during which much of this work was completed.

\newpage
\appendix
\section{Computation of Green's function\label{apA}}

We are looking for Green's function for the following linear BVP
\begin{alignat}{2}\label{AppBVP}
\mathcal{L}v &= \left[1+\Delta R+\Delta \partial_{xx}+\Delta \nu \partial_{xxxx}\right] v=r &\quad v(-1) &= v(1)=v_{xx}(-1)=v_{xx}(1)=0
\end{alignat}
where $r$ will be set by all terms treated explicitly in the time discretization. With these boundary conditions, $\mathcal{L}$ is
symmetric, and has the following spectrum
\begin{alignat}{2}\label{spectrum}
w^{\rm o}_k &=\sin(k \pi x) &\quad \lambda^{\rm o}_k &= 1+\Delta R -\Delta \pi^2 k^2+\Delta \nu \pi^4 k^4\\
w^{\rm e}_k &=\cos\left(\left[k-\frac{1}{2}\right] \pi x\right) &\quad \lambda^{\rm e}_k &= 1+\Delta R -\Delta\pi^2 \left(k-\frac{1}{2}\right)^2+\Delta\pi^4  \nu\left(k-\frac{1}{2}\right)^4
\end{alignat}
where the superscripts denote odd and even. We use the eigenfunction decomposition of Green's function, given by
\begin{equation}\label{Geigen}
G(x,y)=\sum\limits_{k=1}^{\infty} \left\{ \frac{w^{\rm o}_k(x)w^{\rm o}_k(y)}{\lambda^{\rm o}_k} + \frac{w^{\rm e}_k(x)w^{\rm e}_k(y)}{\lambda^{\rm e}_k}\right\} 
\end{equation}
The summations over the odd and even contributions proceeds in a similar fashion, so we will focus on the former. After factorising the denominator as
\begin{equation}\label{roots}
\lambda^{\rm o}_k=\Delta \nu(\pi^2 k^2 - p_+^2)(\pi^2 k^2 - p_-^2);\ p^2_{\pm}=\frac{1}{2\nu}\pm \frac{1}{2\Delta\nu}\sqrt{\Delta^2-4\Delta\nu (1+\Delta R)}
\end{equation} 
we can expand the summation as
\begin{multline}
\sum\limits_{k=1}^{\infty} \frac{w^{\rm o}_k(x)w^{\rm o}_k(y)}{\lambda^{\rm o}_k}=\frac{1}{2\Delta\nu (p^2_+-p^2_-)} \sum\limits_{k=1}^{\infty}  \left\{\frac{\cos(k\pi\, [x-y])}{\pi^2 k^2-p^2_+}
-\frac{\cos(k\pi\, [x-y])}{\pi^2 k^2-p^2_-}\ \ \ \ \ \ \ \right. \\
\shoveright{\left. -\frac{\cos(k\pi\, [x+y])}{\pi^2 k^2-p^2_+}+\frac{\cos(k\pi\, [x+y])}{\pi^2 k^2-p^2_-}\right\}=\ldots=}\\
\frac{1}{2} \Re\left(\frac{1}{2\Delta\nu (p^2_+-p^2_-)} \left\{ \widecheck{f}_{+}(x-y) -\widecheck{f}_{-}(x-y)-\widecheck{f}_{+}(x+y)+\widecheck{f}_{-}(x+y)\right\}\right)
\end{multline}
where the ellipsis corresponds to some tedious manipulations of the sums to bring them into the form of the elementary
inverse semi discrete Fourier transform
\begin{equation}
\widecheck{f}_{\pm}(x)=\left(\frac{1}{\pi^2 k^2-p^2_{\pm}}\right)\raisebox{10pt}{$\!\!\widecheck{\mbox{}}$}=\sum\limits_{k=-\infty}^{\infty}\left(\frac{1}{\pi^2 k^2-p^2_{\pm}}\right)  e^{i\pi k x}=-\frac{\cos( p_{\pm} - p_{\pm}|x|)}{ p_{\pm}\sin( p_{\pm})} 
\end{equation}
As can be seen from this expression, Green's function will exhibit global oscillations if the discriminant in Eq. \ref{roots} is positive so that 
at least one of $p^2_{\pm}$ is real-valued. In that case, the resulting time-stepping scheme will be inaccurate and often unstable, as 
demonstrated in Sec. \ref{tests}.
We will therefore assume that $\Im(p_{\pm})\ne 0$, which means that we impose an upper bound on the time step. 

Combining the odd and even contributions, we obtain Green's function in the compact form
\begin{equation}\label{Green1}
G(x,y)=\frac{1}{\Delta \nu \Im(p^2_+-p^2_-)} 
\Im\left(\frac{-\cos(2 p_+- p_+|x-y|)+\cos( p_+ (x+y))}{p_+ \sin(2 p_+)}\right)
\end{equation} 
In this form, Green's function includes terms as large as $\exp(2\Im(p_+))$ near the boundary, and $p_+$, in turn, grows as $1/\sqrt{\nu}$.
This causes large cancellation errors near the boundaries. A more suitable form is
\begin{multline}\label{Green2}
G(x,y)= Q \left(e^{- b|x-y|}\sin(2 a- a |x-y|-\phi)
-e^{-4 b+ b |x-y|} \sin(2 a - a |x-y|+\phi) \right.\\
\left. -e^{-2 b+ b (x+y)}\sin( a (x+y)-\phi)+e^{-2 b- b (x+y)}\sin( a (x+y)+\phi)\right)
\end{multline}
where we have introduced the auxiliary parameters
\begin{alignat}{2}\label{auxiliaries}
S&=\frac{\Delta}{4\nu [1+\Delta R]} & &\qquad\text{($0<S<1$)} \nonumber\\
a&=\Re(p_+)=\frac{S^{-1/4}}{2\sqrt{\nu}}\sqrt{1+\sqrt{S}} & & \nonumber\\
b&=\Im(p_+)=\frac{S^{-1/4}}{2\sqrt{\nu}}\sqrt{1-\sqrt{S}} & & \nonumber\\
\theta&=\mbox{Arg}(p_+)=\arctan\sqrt{\frac{1-\sqrt{S}}{1+\sqrt{S}}} & &\qquad \text{($0<\theta<\pi/4$)} \nonumber\\
\phi&=\mbox{Arg}\left(\frac{1}{p_+\sin(2 p_+)}\right)=2a-\theta-\frac{\pi}{2}+\arctan\left(\frac{\sin(4a)}{e^{4b}-\cos(4a)}\right) & & \nonumber\\
Q&=\frac{\sqrt{\nu}}{\Delta} \sqrt{\frac{2 S\sqrt{S}}{1-S}} \frac{1}{\sqrt{1-2e^{-4 b}\cos(4 a)+e^{-8 b}}} & &
\end{alignat}
Where $p_+$ has been chosen to lie in the first quadrant in the complex plane.

In this form, Green's function no longer has exponentially large factors. It is immediately clear that, if the ratio between
$\nu$ and $\Delta$ is fixed, then the amplitude of $G$ grows only as $1/\sqrt{\nu}$ for small viscosity. However, there are still
terms of $\mbox{O}(1)$ that cancel near the boundary. To avoid this, we rewrite Green's function near $x=1$ as
\begin{multline}
G(x,y)=Q \sin(a (x-1)) \left[\cos(a (y+1)+\phi) (e^{ b (x-y-4)}+e^{- b (x+y+2)})-\right]\\
   \shoveright{\left.\cos(a (y+1)-\phi) (e^{- b (x-y)}+e^{ b (x+y-2)}) \right]}\\
   \shoveleft{-2 Q \sinh(b (x-1)) \cos(a (x-1)) \left[\sin(a (y+1)-\phi) e^{b (y-1)}+\right.}\\
   \left.\sin(a (y+1)+\phi) e^{- b (y+3)}\right]
\end{multline}
and use the latter form for numerical evaluation if $1-x<1/b$. 

Green's function satisfies
\begin{alignat}{2}\label{map}
G(x,y)&=G(y,x) &\qquad G(-x,-y)&=G(x,y)
\end{alignat}
the latter property being a consequence of an $S_2$ symmetry of BVP \rf{AppBVP}, namely
\begin{alignat*}{2}
( v,\, x,\, t) \ \rightarrow\ (-v,\,-x,\,t)
\end{alignat*}
The BVP is equivariant under this reflection if, and only if, the original boundary conditions are, i.e. if $r=-l$
in \rf{KSIBVP}-\rf{BC1}.
As a consequence of \rf{map}, we only have
to evaluate Green's function in on the domain $0\leq x\leq 1$, $-x\leq y\leq x$. Similar expressions are
readily derived for $G_y(x,y)$ near the boundaries $x=1$ and $y=-1$. 
\bibliographystyle{siam}
\bibliography{KS_Green}
\end{document}